\documentclass[12pt]{article}         
\usepackage{amsfonts}  
\title{First eigenvalue and Maximum principle for fully nonlinear singular operators.}
\author{I. Birindelli, F. Demengel}
\date{}

\newtheorem{theo}{Theorem}[section]
\newtheorem{prop}[theo]{Proposition}
\newtheorem{rema}[theo]{Remark}
\newtheorem{defi}[theo]{Definition}
\newtheorem{cor}[theo]{Corollary}
\newtheorem{lemme}[theo]{Lemma}

\newcommand{\R}{{\rm I}\!{\rm  R}}
\def\R{{\rm I}\!{\rm  R}}
\def\la1{\lambda_1}
\def\grad{\nabla}
\newcommand{\N}{{\bf N}}

\begin{document}
\maketitle
\section{Introduction}

Let $Lu=tr(A(x)D^2u)$ where $A(x)$ is a positive definite matrix satisfying
$mI\leq A(x)\leq MI$ for some positive constants $m$ and $M$. If $\Omega$ is a 
bounded smooth domain of $\R^N$, it is well known that there exists $\bar \lambda$ 
such that:

\begin{itemize}
\item There exists a positive function $\phi$ satisfying
$$\left\{\begin{array}{lc}
L\phi +\bar \lambda\phi=0 & {\rm in}\ \Omega\\
\phi=0 &  {\rm on}\ \partial \Omega.
\end{array}
\right.
$$ 
Furthermore $\bar\lambda$ is the smallest eigenvalue of $-L$ and hence :

\item For any $\lambda<\bar \lambda$ and for any $f\in L^N(\Omega)$ there  
exists a unique $u$ such that
$$\left\{\begin{array}{lc}
Lu + \lambda u=f & {\rm in}\ \Omega\\
u=0 & {\rm on}\ \partial \Omega.
\end{array}
\right.
$$ 
(see \cite{GT}). 
\end{itemize}

Let us recall that $L+\lambda $ {\em satisfies the maximum principle in }
$\Omega$ if any solution of $Lu +\bar \lambda u\leq 0$  in $ \Omega$ which is 
positive on the boundary of $\Omega$ is positive in $\Omega$.

The first eigenvalue of $-L$ in $\Omega$ is characterized by the fact that 
it is the 
supremum of the real numbers  $\lambda$ such that $L+\lambda $ satisfies the maximum
 principle in $\Omega$. (See Protter and  Weinberger  when $L$ is the Laplacian
 and Berestycki, Nirenberg and Varadhan \cite{BNV} for general second order 
operators in general domains.)

\bigskip
On the other hand if $Lu:=\Delta_p u:={\rm div}(|\grad u|^{p-2}\grad u)$ the value
$$\bar\lambda=\inf \frac{\int_\Omega |\grad u|^p}{\int_\Omega | u|^p}$$
has been called the first eigenvalue for $-\Delta_p$ even though strictly speaking it is not (see e.g. \cite{A,L}).
All the same $\bar\lambda$ has the required properties of the eigenvalue:

\begin{itemize}
\item There exists a positive function $\phi$ satisfying
$$\left\{\begin{array}{lc}
\Delta_p\phi +\bar \lambda\phi^{p-1}=0 & {\rm in}\ \Omega\\
\phi=0 &  {\rm on}\ \partial \Omega.
\end{array}
\right.
$$

\item For any $\lambda<\bar \lambda$ and for any $f\in L^p(\Omega)$ there  
exists a unique $u$ such that
$$\left\{\begin{array}{lc}
\Delta_p u + \lambda |u|^{p-2}u=f & {\rm in}\ \Omega\\
u=0 & {\rm on}\ \partial \Omega.
\end{array}
\right.
$$ 
\end{itemize}
It is important to remark  that the higher order term and the zero order term are  homogeneous
of the same degree and that the definition of $\bar \lambda$ is related to the variational nature of
the $p$-Laplacian.

\bigskip

In this paper we introduce a notion of {\em first eigenvalue} for fully nonlinear operators which are non variational but 
homogeneous. Following Berestycki, Nirenberg and Varadhan \cite{BNV} this so called eigenvalue will be defined through the 
{\em maximum principle} but it will have the feature of the ``eigenvalue'' of $-\Delta_p$.

 Indeed we  consider fully nonlinear elliptic 
operators $Lu:=F(\grad u,D^2u)$  which may be singular or degenerate (as the $p$-Laplacian) that satisfy

\begin{itemize}
\item[(F1)] $F(tp,\mu X)=|t|^{\alpha}\mu F(p,X)$, $\forall t\in \R$,   $\mu\in\R^+$, $\alpha>-1$

\item[(F2)] $a |p|^\alpha {\rm tr} N\leq F(p,M+N)-F(p,M)\leq A |p|^\alpha{\rm tr} N$ for $0<a\leq A$, $\alpha>-1$ and $N\geq 0$.
\end{itemize}
The class of operators satisfying (F1) and (F2) is large and includes
$$ F(\grad u,D^2u)=|\grad u|^\alpha {\cal M}_{a,A}(D^2u)$$
where $\alpha>-1$ and ${\cal M}_{a,A}$ is one of the Pucci operators.
$$F(\grad u,D^2u)=\Delta_p u$$
with $\alpha=p+1$.
In \cite{BD} many examples of operators satisfying (H2) are given.

Of course the right notion of solution in this context will be 
that of viscosity solution (see \cite{CIL}) suitably adapted to our contest. Let us remark that even 
if $u$ is $C^2$, when $\alpha<0$, $F(\grad u,D^2u)$ is not define for $\grad u=0$.

\bigskip 

Before going into details, let us mention that a certain number of interesting
 papers have appeared that treat viscosity solutions for equations involving the 
$p$-Laplacian. In  fact in \cite{JLM,JM} Juutinen, Lindqvist and Manfredi opened the 
way to this topic. We would like to emphasize that in those 
papers the point of view is really on the $p$-Laplacian and the variational 
structure of the $p$-Laplacian is used. 
This is not the case here since the operators we consider are fully nonlinear.

\bigskip

The first key ingredient is the following:

\begin{theo}\label{maxp}
Suppose that $\Omega $ is a bounded open piecewise ${\cal C}^1$ domain of $\R^N$.
Suppose that for $\lambda\in\R$ there exists a function $v>0$ such that $\ F(\grad v,D^2v)+\lambda v^{\alpha+1}\leq 0$ in $\Omega$. 
Then, for $\tau<\lambda$,  every viscosity solution of 

$$\left\{\begin{array}{c}
F(\grad \sigma,D^2 \sigma)+\tau\vert \sigma\vert^{\alpha} \sigma \geq 0\ {\rm in }\ \Omega\\
 \sigma \leq 0\ {\rm on}\ \partial \Omega
\end{array}
\right.
$$
satisfies 
$\sigma\leq 0$ in $\Omega$.
\end{theo}
When $\tau<0$ this result was obtained in \cite{BD} for the operators considered here, of course for a large 
class of elliptic operators see \cite{CIL} and the references therein. It is well known that to prove maximum principles or comparison principles for viscosity solutions one needs to double the variables and consider
the function $\psi(x,y)=u(x)-v(y) +\phi(x,y)$ where $\phi$ is an
appropriate $C^2$ function (see \cite{CIL}). On the other hand here
instead of considering a difference of sub and super solutions we
consider the ratio of $\sigma$ and $v$.

This theorem allows us to define
$$\bar\lambda=\sup\{ \lambda\in \R,\ \ \exists \  \phi>0\  \ {\rm in}\ \Omega,\ F(\grad \phi,D^2\phi)+\lambda \phi^{\alpha+1}\leq 0 \mbox{ in the viscosity
sense}
\ \}.$$
In other words if we denote by  $I_\alpha(u)=|u|^\alpha u$, $\bar\lambda$ is the supremum of the value $\lambda$ such that 
$F+\lambda I_\alpha$ satisfies the maximum principle in $\Omega$.
The main aim of this paper is to convince the reader that it is {\em correct} to call $\bar\lambda$  the first eigenvalue of $-F$
 in $\Omega$.

\bigskip
Clearly the set 
$$E=\{ \lambda\in \R,\ \ \exists \  \phi>0\  \ {\rm in}\ \Omega,\ F(\grad \phi,D^2\phi)+\lambda \phi^{\alpha+1}\leq 0 \mbox{ in the viscosity
sense}
\ \}$$ 
is an interval;  in fact it is an interval which is bounded from above since

\begin{prop}
Suppose that $R$ is the radius of the largest ball contained in the bounded set $\Omega$. Then, there exists some constant $C$ which depends only on $N$ and $\alpha$, such that 
$$\bar\lambda\leq {C\over R^{\alpha+2}}.$$
\end{prop}
The value $\bar\lambda$ has the following features that justify the name of eigenvalue:

\begin{theo}\label{exiev}
There exists $\phi$ a continuous  positive viscosity solution of 
$$\left\{\begin{array}{lc}
F(\grad \phi,D^2\phi) +\bar\lambda\phi^{\alpha+1}=0 & {\rm in}\ \Omega\\
\phi=0 &  {\rm on}\ \partial \Omega.
\end{array}
\right.
$$ 
\end{theo}
Furthermore
\begin{theo}\label{exif}
For $\lambda<\bar\lambda$ if $f<0$ in $\Omega$ and bounded then there
exists a unique $u$ nonnegative viscosity solution of
\begin{equation}\label{eqif}\left\{\begin{array}{lc}
F(\grad u,D^2u) + \lambda u^{\alpha+1}=f & {\rm in}\ \Omega\\
u=0 &  {\rm on}\ \partial \Omega.
\end{array}
\right.
\end{equation}        
\end{theo}
Let us notice that in order to prove the Theorems \ref{exiev} and \ref{exif} we need to obtain some estimates which are interesting in their own right:
\begin{theo} \label{hol}Suppose that $f$ is a bounded function in $\bar\Omega$ then if $u$ is a bounded nonnegative
viscosity solution of 
$F(\grad u,D^2u)=f$ in $\Omega$, it is H{\"o}lder continuous:
$$|u(x)-u(y)|\leq M|x-y|^\gamma.$$
\end{theo}

In \cite{IL}, among other results, Ishii and Lions prove H{\"o}lder and Lipschitz estimates for a class of second order elliptic operators that do not include the operators considered here but the proof of Theorem \ref{hol}
is inspired by \cite{IL}. We also obtain local Lipschitz regularity using the H{\"o}lder regularity of the solution.
Recently Arizawa and Capuzzo Dolcetta in  \cite{ACD} have obtained  H{\"o}lder's regularity for solutions of other
degenerate elliptic fully nonlinear operators.

In the case $\alpha=0$ H{\"o}lder's regularity is proved by Caffarelli Cabr{\'e} with a different proof that requires Harnack and Alexandrov-Bakelman-Pucci estimates.

\bigskip
After the completion of this work we learned that the eigenvalue problem for 
Pucci's operators has already been treated; this would correspond to the case $\alpha=0$ here. An initial work concerning the radial case
was completed by Felmer and Quaas \cite{FQ} then Quaas treated the case of  general domains \cite{Q}. Later the work was completed by  
Busca, Esteban and Quaas in \cite{BEQ}. In this interesting paper they denote by $\mu_1^+$ the eigenvalue here denoted $\bar \lambda$, but they also define 
$$\mu_1^-=\sup\{\mu,\ \exists \psi<0 \ {\rm in}\ \Omega : {\cal M}^+_{a,A}(\psi)+\mu\psi\geq 0\}.$$
This could be done in our case as well but we have limited ourselves to positive solutions. On the other hand, a priori estimates for ${\cal M}^+_{a,A}$ have been given in \cite{CC} that allow them to take solution in $W^{2,n}$. These estimates are not known for singular operators.
They also consider interesting bifurcation problems.

\bigskip
Let us mention some open problems:

\begin{itemize}
\item {\bf Simplicity of the eigenfunction.} The first eigenfunction $\phi$ is simple for linear second order elliptic operators, for the Pucci's operators  and for the $p$-Laplacian,
it would be interesting to know if this is true also in the case treated here i.e. suppose that $\psi>0$ is another eigenfunction
does this imply that there exists $t\in\R^+$ such that $\psi=t\phi$?

\item {\bf Fredholm alternative} By the definition of $\bar\lambda$ if $f<0$ in $\bar\Omega$ then there are no positive solutions of
$$\left\{\begin{array}{lc}
F(\grad u,D^2u) +\bar\lambda |u|^{\alpha}+1 =f & {\rm in}\ \Omega\\
u=0 &  {\rm on}\ \partial \Omega.
\end{array}
\right.
$$ 
Is it still true if $f\leq 0$? are there solutions that change sign?

\item {\bf $\bar\lambda$ is isolated.} Suppose that $\bar\lambda+\varepsilon \lambda>\bar \lambda$.
Is it possible to prove that  there exists a solution of (\ref{eqif}) if $\varepsilon$ is sufficiently small?

\end{itemize}

In the next section we state the precise hypothesis on the fully non linear operator $F$ and we give 
the notion of ``viscosity solution'' adapted to the operators considered here. 
In the third section we prove the maximum principle (Theorem \ref{maxp}) and a comparison principle. In the fourth section we give global H{\"o}lder and local 
Lipschitz estimates for the solutions. 
Finally in the last section we prove different existence results including that of a first eigenfunction.
Some properties of the distance function are proved in the appendix.

{\bf Acknowledgment:}  This work was completed while I. Birindelli was
visiting the University of Cergy-Pontoise with the  support of the
Laboratoire d'Analyse, 
g{\'e}om{\'e}trie et Mod{\'e}lisation and of CNR- Short term program 2004, she would like to thank both institutions.
Both authors wish to thank I. Capuzzo Dolcetta for interesting conversations.

\section{Preliminaries}

Let $\alpha$ be some real number, $\alpha>-1$, and let $F$ be a
 fully nonlinear and singular or degenerate operator 
$$F(\grad u,D^2 u) $$
where $F: \R^N-\{0\}\times S$,  $S$ is the set of symmetric matrix in
$\R^N$, and  we consider the following hypothesis

\bigskip

(H1) $F(tp,\mu X)=|t|^{\alpha}\mu F(p,X)$, $\forall t\in \R$, $\mu\in\R^+$, $\alpha>-1$ 
and $F(p,X)\leq F(p,Y)$ for any $p\neq 0$, and $X\leq Y$.

\bigskip

In most of the paper the operator satisfies 
also the following  hypothesis

\bigskip

(H2) $a |p|^\alpha {\rm tr} N\leq F(p,M+N)-F(p,M)\leq A |p|^\alpha{\rm tr} N$ for $0<a\leq A$, $\alpha>-1$ and $N\geq 0$.

\bigskip
\noindent When only (H1) is required it will be stated explicitly.

\noindent Let us recall that (H2) implies
$$|p|^\alpha{\cal M}^+_{a,A}(M)\geq F(p,M)\geq |p|^\alpha{\cal M}^-_{a,A}(M),$$
where, if $e_i$ are the eigenvalues of $M$

$${\cal M}^+_{a,A}(M)=a \sum_{e_i<0} e_i+A\sum_{e_i>0} e_i \quad {\rm and}$$ 

$${\cal M}^-_{a,A}(M)=A \sum_{e_i<0} e_i+a\sum_{e_i>0} e_i $$
are the Pucci operators (see e.g. \cite{CC}).

\begin{rema}
Let us observe that if $F$ satisfies (H2), $G(p,X) = -F(p, -X)$ satisfies
(H2). With this remark, defining 
$$\bar\lambda^- = \sup\{ \mu, \exists \ \phi<0\ , F(\nabla \phi, D^2\phi)+
\mu |\phi|^\alpha \phi\geq 0\}$$
and observing that $\bar \lambda ^-= \bar\lambda (G)$ one gets symmetrical
results to those enclosed in the sequel for the value $\bar\lambda$.  
\end{rema}

\bigskip

We need first to extend the definitions employed in \cite{BD}. Let us recall 
first the definition of viscosity continuous sub or
super solutions for operators that satisfy $(H2)$ and hence may be singular when $\grad u=0$.

It is well known that in dealing with viscosity respectively  sub and super
solutions one works with  
$$u^\star (x) = \limsup_{y, |y-x|\leq r} u(y)$$

and $$u_\star (x) = \liminf_{y, |y-x|\leq r} u(y).$$
 It is easy to see that 
$u_\star \leq u\leq u^\star$ and $u^\star$ is uppersemicontinuous (USC)
$u_\star$ is lowersemicontinuous (LSC).  See e.g. \cite{CIL, I}.

\begin{defi}\label{def1}

 Let $\Omega$ be an open set in
$\R^N$, then
$v$  bounded on $\overline{\Omega}$ is called a viscosity super-solution
of
$F(\grad v,D^2v)=g(x,v)$ if for all $x_0\in \Omega$, 

-Either there exists an open ball $B(x_0,\delta)$, $\delta>0$  in $\Omega$
on which $v$ is constant and equal to $c$ then
$g(x, c)\geq  0$

-Or
 $\forall \varphi\in {\cal C}^2(\Omega)$, such that
$v_\star-\varphi$ has a local minimum on $x_0$ and $\nabla \varphi(x_0)\neq
0$, one has
\begin{equation}F( \nabla \varphi(x_0),
 D^2\varphi(x_0))\leq g(x_0, v_\star(x_0)).
\label{eq1}\end{equation}

\bigskip
\noindent Of course $u$ is a viscosity sub-solution if
for all $x_0\in \Omega$,

-Either there exists a ball $B(x_0, \delta)$, $\delta>0$ on which  $u$ is constant and equal to $c$
then
$g(x, c)\leq 0$,

-Or  $\forall
\varphi\in {\cal C}^2(\Omega)$, such that
$u^\star-\varphi$ has a strict local maximum on $x_0$ and
$\nabla\varphi(x_0)\neq 0$, one has
\begin{equation}F( \nabla \varphi(x_0),
D^2\varphi(x_0))\geq g(x_0, u^\star(x_0)). \label{eq2}\end{equation}
\end{defi}

See e.g. \cite{CGG} and \cite{ES} for similar definition of viscosity solution for 
equations with singular operators.

For convenience we recall  the definition of semi-jets given e.g. in 
\cite{CIL}

\begin{eqnarray*}
J^{2,+}u(\bar x) &= &\{ (p, X)\in \R^N\times S, \ u(x)\leq
u(\bar x)+
\langle p, x-\bar x\rangle + \\
&+& {1\over 2} \langle X(x-\bar x), x-\bar
x)\rangle
+ o(|x-\bar x|^2) \}
\end{eqnarray*}
and 

\begin{eqnarray*}J^{2,-}u(\bar x)& =& \{ (p, X)\in \R^N\times S, \
u(x)\geq u(\bar x)+
\langle p, x-\bar x\rangle +\\
&+& {1\over 2} \langle X(x-\bar x), x-\bar
x)\rangle
+ o(|x-\bar x|^2\}.
\end{eqnarray*}
In the definition of viscosity solutions the test functions can be
substituted by the elements of the semi-jets in the sense that if $(p,X)\in J^{2,-}u(\bar x)$ then 
$\phi(x)=u(\bar x)+
\langle p, x-\bar x\rangle 
+ {1\over 2} \langle X(x-\bar x), x-\bar
x)\rangle$ is a test function for a super solution $u$ at $\bar x$ and similarly if 
$(q,Y)\in J^{2,+}u(\bar x)$ then 
$\varphi(x)=u(\bar x)+
\langle q, x-\bar x\rangle 
+ {1\over 2} \langle Y(x-\bar x), x-\bar
x)\rangle$  is a test function for a sub solution $u$ at $\bar x$.

\section{Maximum principle and comparison results}

As we pointed out in the introduction, we want to generalize the concept of 
eigenvalue for the  Dirichlet problem in a bounded domain $\Omega$ associated to the operator $L(u) = F(\nabla u, D^2 u)$ satisfying (H2) .
 It will be defined following the main ideas introduced in 
\cite{BNV} for linear uniformly elliptic operators.

\bigskip
 In the introduction we have defined the eigenvalue $\bar \lambda$.
In view of the definition given in the previous section the correct definition of the set $E$ becomes 
$$E = \{ \lambda\in \R,\ \ \exists \  \phi, \ \phi_\star>0  \ {\rm in}\
\Omega,\ F(\grad \phi,D^2\phi)+\lambda \phi^{\alpha+1}\leq 0 \mbox{ in the viscosity
sense}
\ \}.$$
Throughout the paper we shall denote 
$$\bar\lambda = \sup E.$$

\bigskip

\begin{rema}
Of course $E$ is non empty, since $0$ obviously belongs to $E$. Moreover, if $\lambda\in E,$ every $\lambda^\prime < \lambda $ is also in $E$.
\end{rema}
The next proposition proves that $\bar\lambda\in \R^+$. In the last
section we shall prove that $\bar\lambda$ plays the role  of the first
eigenvalue.

\begin{prop}\label{pro1}
Suppose that $R$ is the radius of the largest ball contained in  $\Omega$ and suppose that $F$ satisfies (H2). Then, there exists some constant $C$ which depends only on $a$, $A$, $N$ and $\alpha$, such that 
$$\bar\lambda\leq {C\over R^{\alpha+2}}.$$
\end{prop}
This proposition is a consequence of the maximum principle stated in the following Theorem and Lemma \ref{lem1} below.

\begin{theo}\label{maxp}
Suppose that $\Omega $ is a bounded open
 piecewise ${\cal C}^1$ domain of $\R^N$.
Suppose that $\tau< \bar\lambda$ and that $F$ satisfies (H1), 
then every viscosity solution of 
$$
\left\{\begin{array}{lc}
F(\grad \sigma,D^2\sigma)+\tau\vert \sigma\vert^{\alpha} \sigma \geq 0 & {\rm in }\ \Omega\\
 \sigma \leq 0 & {\rm on}\ \partial \Omega
\end{array}
\right.
$$
satisfies 
$\sigma\leq 0$ in $\Omega$.
\end{theo}
{\bf Remark:} In this theorem we don't require $F$ to satisfy (H2), but 
only (H1).

An immediate consequence of Theorem \ref{maxp} is
\begin{cor} Suppose that  $F$ satisfy (H2).
If $\lambda<\bar\lambda$ and $F(-p,-X)=-F(p,X)$ then every solution of 
\begin{equation}\label{cor34}\left\{\begin{array}{lc}
F(\nabla \phi, D^2\phi)+\lambda |\phi|^\alpha \phi = 0 & {\rm in }\ \Omega\\
 \phi = 0 & {\rm on}\ \partial \Omega
\end{array}
\right.
\end{equation}
which is zero on the boundary,  is identically zero. 

Let 
$$
\lambda^+= \{ \lambda\in \R,\ \ \exists \  \phi, \ \phi_\star>0  \ {\rm in}\
\Omega,\ |\grad \phi|^\alpha{\cal M}_{a,A}^+(D^2\phi)+\lambda \phi^{\alpha+1}\leq 0 \mbox{ in the viscosity
sense}
\ \},
$$ then for any $\lambda<\lambda^+$ the solution of (\ref{cor34}) is zero.

\end{cor}

{\em Proof:} In the first case, both $\phi$ and $-\phi$ are solutions of the equation and
this implies that they are both negative. 
In the general case the hypothesis on $F$ implies that $\phi$ is a subsolution of 
\begin{equation}\label{eqif2}
|\grad \phi|^\alpha{\cal M}^+_{a,A}(D^2\phi) +
\lambda|\phi|^{\alpha}\phi\geq 0\ {\rm in}\ \Omega,
\end{equation}
and hence $\phi\leq 0$. On the other hand $\phi$  is a supersolution of
$$|\grad \phi|^\alpha{\cal M}^-_{a,A}(D^2\phi) +
\lambda|\phi|^{\alpha}\phi\leq 0\ {\rm in}\ \Omega,$$ and therefore
$-\phi$ is a supersolution (\ref{eqif2}) and $-\phi$ is also negative.
This conclude the proof.

\begin{lemme}\label{lem1}
Suppose that $\Omega= B(0,R)$, and let  $q= {\alpha+2\over \alpha+1}$ and 
$$\sigma = {1\over 2q} (\vert x\vert^q-R^q)^2.$$
Let $F$ satisfy (H2).
Then there exists some constant $C$ which depends only on $a$, $A$, $N$ and $\alpha$ such that  
$$\sup_{x\in B(0,R)} {-F(\grad \sigma,D^2\sigma)\over \sigma^{\alpha+1}}\leq {C\over R^{\alpha+2}}$$
\end{lemme}
{\it Proof of Proposition \ref{pro1}.}
Suppose that Theorem \ref{maxp} and Lemma \ref{lem1}  hold. Without loss of generality 
we can suppose that $B(0,R)\subset\Omega$. We shall prove that 

$$\bar \lambda \leq 
\sup_{x\in B(0,R)} {-F(\grad \sigma,D^2\sigma)\over \sigma^{\alpha+1}}= \tau,$$
by Lemma \ref{lem1} this ends the proof.

Suppose by contradiction that $\tau< \bar\lambda  $ and 
let $u=\sigma$ for $|x|\leq R$ and 0 elsewhere. Then one 
 would have 
$$F(\grad u,D^2u) +\tau \vert u\vert ^\alpha u \geq 0\  {\rm in} \ \Omega.$$

Indeed, for $|x|\leq R$, $u$ is a solution  by the definition of $\tau$, for $|x|>R$ 
the definition of viscosity solution gives the result  immediately and for $|x|=R$ all the test  functions have zero gradient and so they don't need to be tested.
Now since $u = 0$ on $ \partial \Omega$, this would imply by Theorem 
\ref{maxp} that 
 $u \leq 0$ in $\Omega$, a contradiction with the definition of $\sigma$ which is positive inside 
the ball. 
This ends the proof of Proposition \ref{pro1}.

\bigskip
\noindent
{\it Proof of Lemma \ref{lem1} :}

Let $g(r) = \sigma (|x|)$. 
The computation of $g^\prime (r)$ gives 
$g^\prime (r) = r^{2q-1}-r^{q-1} R^q$ and
$$g^{\prime\prime }(r) = (2q-1) r^{2q-2}-(q-1) r^{q-2} R^q.$$
Clearly $g'\leq 0$ while $g''\leq 0$ for $r\leq \left(\frac{q-1}{2q-1}\right)^{1\over q}$
and positive elsewhere. Hence by condition (H2) and using the fact that for radial functions the
eigenvalues of the Hessian are $\displaystyle{\frac{g'}{r}}$ with multiplicity N-1 and $g''$ (see \cite{CL}),  
$$F(\grad \sigma,D^2\sigma) \leq |g'|^\alpha{\cal M}^+_{a,A}=|g'|^\alpha\left[ag^{\prime\prime}(r)+ a({N-1\over r}) g^\prime(r)\right] $$ or
$$F(\grad \sigma,D^2\sigma) \leq |g'|^\alpha{\cal M}^+_{a,A}=|g'|^\alpha\left[A g^{\prime\prime}(r)+ a({N-1\over r}) g^\prime(r)\right] .$$
In both cases

$$F(\grad \sigma,D^2\sigma) \leq |g'|^\alpha r^{q-2} (B_1r^q-B_2R^q)$$
with either
$B_1 = a(N+2q-2)$ and $B_2 = a(N+q-2)$ or

 $B_1=A(2q-1)+a(n-1)$ and $B_2=A(q-1)+a(N-1)$.
Hence one gets:
$$\frac{F(\grad \sigma,D^2\sigma)}{\sigma^{\alpha+1}}\leq
%-{|g^\prime (r)|^\alpha (g^{\prime\prime}+ ({N-1\over r} g^\prime  )(r)\over (g(r)^{\alpha+1}} = 
- {r^{q(\alpha+1)-\alpha-2} (-B_1 r^q+B_2 R^q)\over (R^q-r^q)^{\alpha+2}}.$$
Let
$$\varphi(r) = { (-B_1 r^q+B_2 R^q)\over (R^q-r^q)^{\alpha+2}},$$
since $q = {\alpha+2\over \alpha+1}$  one has that 
$$\frac{F(\grad \sigma,D^2\sigma)}{\sigma^{\alpha+1}}\leq \sup \varphi(r).$$
It is easy to see that $\sup \varphi(r)=\frac{C}{R^{\alpha+2}}$.
% since
%$$\varphi(0) = {B_2\over R^{\alpha+2}}$$
%and $\lim_{r\rightarrow R} {\varphi} (r) = -\infty$
%and $\varphi^\prime (r) = r^{q-1} {q(B_1+B_2(\alpha+2)) R^q-q(\alpha+3) r^q\over (R^q-r^q)^{\alpha+3}}$
%which is zero for 
%$$r_0 = {B_1+(\alpha+2) B_2\over \alpha+3} R$$
%and then 
%$$\sup \varphi \leq \varphi (r_0)= {C\over R^{\alpha+2}}.$$
 
 This ends the proof of Lemma \ref{lem1}.

\bigskip

\noindent{\bf Proof of Theorem \ref{maxp}}. We assume that $\tau<\bar\lambda$. Then taking $\lambda $ such that $\tau<\lambda<\bar\lambda$,
there exists  $v$, a viscosity sub solution of 
$$F(\grad v,D^2v)+\lambda v^{1+\alpha}\leq 0\ \  {\rm in}\ \Omega,$$
with $v_\star>0$ in $\Omega$. Suppose that $\sigma$ is a 
viscosity solution of 
$$F(\grad \sigma,D^2\sigma)+\tau \vert\sigma\vert^\alpha \sigma \geq 0\ \  {\rm in}\ \Omega,$$
and $\sigma\leq 0$ on $\partial \Omega$. We need to prove that
$\sigma\leq 0$ in $\Omega$. 
It is sufficient  to prove that $\sigma^\star \leq 0$. 
Using the definition of viscosity solutions one can assume without loss
of generality that  
$\sigma
\in USC(\bar\Omega)$ and
$v\in LSC(\bar\Omega)$ and hence drop the stars.

Let us suppose by contradiction that ${\sigma (x)\over v(x)}$ has a positive supremum 
inside $\Omega$.  For some $q>2$ let
us consider the function 
$$\psi_j(x,y)=\frac{\sigma(x)}{v(y)}-{j\over q v(y)}|x-y|^q$$
which is uppersemicontinuous.
Then $\psi_j$ also has a positive supremum   achieved in some couple of
points $(x_j, y_j)\in \Omega^2$. One easily has  that 
$(x_j, y_j)\rightarrow (\bar x, \bar x)$, $\bar
x\in \Omega$ which is a supremum for ${\sigma \over v}$.
 One can also prove  that $j|x_j-y_j|^q\rightarrow 0$, and that $\bar x$ is a continuity point for $\sigma$. For that aim remark
that 
$${\sigma (x_j)-{j\over q} |x_j-y_j|^q\over v(y_j)}\geq {\sigma (\bar x)\over v(\bar x)}$$
and using the lowersemicontinuity of $v$ on $\bar x$ together with $\lim {j\over q} |x_j-y_j|^q=0$
one gets 
$$\liminf\sigma (x_j)\geq \sigma (\bar x).$$

Assume for the moment that $x_j\neq y_j$ for $j$ large enough. Take $j$ large enough in order that 
$$\sigma (x_j)^{1+\alpha} \geq {3\sigma (\bar x)^{1+\alpha}\over 4}$$
and 
$${j\over q} |x_j-y_j|^q\leq {\sigma (\bar x)^{1+\alpha} (\lambda-\tau)\over 4\lambda}.$$
Using $\psi_j(x,y)\leq \psi_j (x_j, y_j)$, 
one gets that 
\begin{equation}\label{eqif3}\sigma (x) v(y_j)-v(y) \left (\sigma (x_j)-{j\over q} |x_j-y_j|^q\right)\leq v(y_j) {j\over q} |x-y|^q.
\end{equation}
We now define
$$\beta_j = \sigma (x_j)-{j\over q} |x_j-y_j|^q$$
and then after some simple calculation (\ref{eqif3}) becomes

\begin{eqnarray}
&&
\left(\sigma(x+x_j) - \sigma (x_j) - j|x_j-y_j|^{q-2} (x_j-y_j. x)\right) v(y_j)+\nonumber\\
& &-\left(v(y+y_j)-v(y_j)  
-j |x_j-y_j|^{q-2} (x_j-y_j. y){v(y_j)\over \beta_j}\right)\beta_j \label{mpf}\\
&\leq&  v(y_j)\left({j\over q} |x_j+x-y_j-y|^q-{j\over q}
|x_j-y_j|^q-j|x_j-y_j|^{q-2} (x_j-y_j, x-y)\right).\nonumber
\end{eqnarray}
We define the functions 
$$U(x) = \left(\sigma(x+x_j)-\sigma (x_j)-j|x_j-y_j|^{q-2} (x_j-y_j. x)\right)v(y_j)
$$
and 
$$V(y) = -\left(v(y+y_j)-v(y_j)-j|x_j-y_j|^{q-2} (x_j-y_j. y){v(y_j)\over \beta_j}\right)\beta_j$$
then (\ref{mpf}) can be written:
$$U(x)+ V(y)\leq (x,y) A (x,y)$$
with 
$$A = j v(y_j)\left(\begin{array}{cc}
D_j & -D_j\\
-D_j& D_j
\end{array}
\right)$$
and 
$$D_j = 2^{q-3} q |x_j-y_j|^{q-2} \left(I+ \frac{(q-2)}{|x_j-y_j|^2} (x_j-y_j)\otimes (x_j-y_j)\right).$$
Then using Theorem 3.2' in \cite{CIL}
one gets that there exist $X_j$ and $Y_j$ such that
$$\left(j|x_j-y_j|^{q-2} (x_j-y_j), {X_j\over v(y_j)}\right)\in J^{2,+} \sigma(x_j)$$
and 
$$\left(j|x_j-y_j|^{q-2} (x_j-y_j){v(y_j)\over \beta_j},{-Y_j\over \beta_j}\right)\in J^{2,-}v(y_j)$$
with, for some $\varepsilon>0$.
$$\left(\begin{array}{cc}
X_j & 0\\
0&Y_j\end{array}
\right)\leq  A+\varepsilon A^2.$$
In particular 
$$X_j+ Y_j\leq 0.$$
We can conclude using the 
 fact that  $v$ and $\sigma$ are respectively a super and  a sub  solution and the properties of $F$. Precisely we have obtained

\begin{eqnarray*}
-\tau \sigma (x_j)^{1+\alpha}&\leq& F(j|x_j-y_j|^{q-2} (x_j-y_j),{X_j\over v(y_j)})\\
&\leq & F(j|x_j-y_j|^{q-2} (x_j-y_j),{-Y_j\over v(y_j)})\\
&\leq & {\beta_j^{1+\alpha}\over v(y_j)^{1+\alpha}} F(j|x_j-y_j|^{q-2} (x_j-y_j){v(y_j)\over \beta_j}, {-Y_j\over \beta_j})\\
&\leq & -\lambda \beta_j^{1+\alpha}=-\lambda[\sigma(x_j)-
\frac{j}{q}|x_j-y_j|^q]^{1+\alpha}.
\end{eqnarray*}
This gives a contradiction, 
indeed by passing to the limit the inequality becomes
$$-\tau\sigma^{\alpha+1}(\bar x)\leq -\lambda \sigma^{\alpha+1}(\bar x).$$

\bigskip

It remains to prove that $x_j\neq y_j$ for $j$ large enough. 
If one assumes that $x_j= y_j$ one has 
$$\sigma (x_j)\geq \sigma (x)-{j\over q} |x_j-x|^q$$
and
$$v(x)\geq v(x_j)-{jv(x_j)|x_j-x|^q\over q \sigma (x_j)}. $$
In that case one uses Lemma 2.2 in \cite {BD} to get a contradiction. 
This ends the proof of Theorem \ref{maxp}.

 \bigskip
Let us recall that in \cite{BD} for $\lambda=0$ we give a  comparison 
principle for continuous viscosity solutions. It is not difficult to see 
that it can be extended to bounded viscosity solutions. We now prove a 
further extension adapted to our context. 

\begin{theo} \label{comp}
Suppose that $\lambda< \bar\lambda$, $f\leq 0$, $f$ is
upper semicontinuous and
$g$ is lower semicontinuous  with $f\leq  g$ and 

- either $f\leq-c<0$ in 
$\Omega$, 

- or  $g(\bar x)>0$ on every point $\bar x$ such that  $f(\bar x)=0$. 

\noindent Suppose that there exist 
$v$   nonnegative   viscosity sub solution of 
$$F(\grad u,D^2u)+\lambda u^{1+\alpha} = f$$
and $\sigma$   nonnegative   viscosity super solution of 
$$F(\grad \sigma,D^2\sigma)+\lambda \sigma^{1+\alpha} \geq g$$
satisfying $\sigma\leq v$ on $\partial \Omega$. 

Then $\sigma \leq v$ in $\Omega$. 
\end{theo}

As a consequence one has 
\begin{cor}
Suppose that $\lambda\leq \bar\lambda$, there exists at most one nonnegative viscosity solution of 
$$\left\{\begin{array}{lc}
F(\grad v,D^2v)+\lambda v^{1+\alpha} = f & {\rm in}\   \Omega\\
v=0  & {\rm on}\  \partial \Omega
\end{array}\right.
$$
  for $f<0$ and continuous.
\end{cor}

{\it Proof of Theorem \ref{comp}}
First using the strict maximum principle (see \cite{BDW} ) one gets
that $F(\grad v,D^2v)\leq 0$
$v\geq 0$, and since $v$ is not identically zero, $v_\star>0$ in
$\Omega$. Without loss of generality one can assume that $\sigma$
and $v$ are respectively USC and LSC. 

Suppose by contradiction that
$\sigma > v$ somewhere in
$\Omega$. The supremum of the function ${\sigma\over v}$  on $\partial \Omega$  is less than $1$, then its supremum
is achieved inside $\Omega$. Let
$\bar x$ be a point such that 
$$1<{\sigma(\bar x)\over v(\bar x)}= \sup_{x\in \overline{\Omega}} {\sigma(x)\over v(x)}.$$

Doing exactly the same construction  as in the proof of Theorem \ref{maxp} we similarly
get : 

\begin{eqnarray*}
g(x_j)-\lambda \sigma (x_j)^{1+\alpha} &\leq& F(j|x_j-y_j|^{q-2} (x_j-y_j), {X_j\over v(y_j)})\\
 &\leq & {\beta_j^{1+\alpha}\over v(y_j)^{1+\alpha}}  F(j|x_j-y_j|^{q-2} (x_j-y_j){v(y_j)\over \beta_j}, {-Y_j\over \beta_j})\\
&\leq& -\lambda \beta_j^{1+\alpha} + {\beta_j^{1+\alpha}\over v(y_j)^{1+\alpha}} f(y_j).
\end{eqnarray*}
Passing to the limit we obtain
\begin{equation}\label{36}g(\bar x)\leq  \left(\frac{\sigma(\bar x)}{v(\bar x)}\right)^{\alpha+1} f(\bar x).
\end{equation}
Either $f(\bar x)=0$ and $g(\bar x)>0.$ but this contradicts (\ref{36}) or
 $f(\bar x)<0$,
and then (\ref{36}) becomes
$$0<f(\bar x)\left[1-\left(\frac{\sigma(\bar x)}{v(\bar x)}\right)^{\alpha+1}\right]\leq f(\bar x)-g(\bar x)\leq 0,$$
also a contradiction.

\noindent This conclude the proof.

\section{H{\"o}lder and Lipschitz  regularity}

In all this section we assume that $F$ satisfies (H2) and $\Omega$ is a $C^2$ bounded domain.

Suppose that $u$ is a  viscosity solution of 
\begin{equation}\label{Heq}
\left\{\begin{array}{lc}
F(\grad u,D^2u) =f & {\rm in}\ \Omega \\
u=0 & {\rm on}\ \partial \Omega.
\end{array}
\right.
\end{equation}

\begin{theo}\label{H1} Let $f$ be a bounded function in $\overline\Omega$.
 Let $u$ be a  non negative viscosity solution of (\ref{Heq}) with $\Omega$
a ${\cal C}^2$ domain.   Then for any $\gamma\in (0,1)$, there exists
$C>0$,
 such that 
$$|u(x)-u(y)|\leq C|x-y|^\gamma.$$
 \end{theo}

Using the result of Theorem \ref{H1}, one obtains the following stronger result :
\begin{theo}\label{L}
 Let $u$ be a non negative viscosity solution of 
equation (\ref{Heq}), then  $u$ is locally  Lipschitz continuous when $f$ is bounded.
\end{theo}

{\em Proof of Theorem \ref{H1}:} The proof relies on ideas
used to prove H{\"o}lder and Lipschitz  estimates in \cite{IL}.
 
First we will prove that $u$ is H{\"o}lder near the boundary  using  the regularity of the 
boundary and of the distance function near the boundary.

\bigskip
Let $d(x)$ be the function $d(x, \partial \Omega)=\inf\{|x-y|,\ {\rm for} \ y\in\partial\Omega\}$. 

\noindent {\bf Claim:} {\em $\exists \delta >0$, there exist $M_o>0$, $1>\gamma>0$ and $\bar r$ such that 
$u(x)\leq M_od(x)^\gamma$  for $d(x)\leq \delta$ }.

In order to prove the claim we need to show that $g(x)=d(x)^\gamma $ is  a super solution of (\ref{Heq})
in $$\Omega_\delta = \{x\in \Omega, d(x, \partial \Omega) < \delta\}.$$

It is well known (see \cite{F,Fo,KP}),  that $d$ is ${\cal C}^2$ on
$\Omega_\delta$ for $\delta$ small enough since $\partial\Omega$ is
${\cal C}^2$.  Furthermore the ${\cal C}^2$ norm of $d$ is bounded. Then for
$\delta$ small enough and $d(x)< \delta$,
$$F(\grad g,D^2g)\leq \gamma ^{1+\alpha}d^{(\gamma(\alpha+1)-\alpha-2)}(\gamma-1+cd(x)|D^2d(x)|_\infty)\leq -\epsilon<0$$
for some constant $c$ which depends on $a$ and $A$ and some 
constant $\varepsilon>0$ which depends on $\gamma, N$, $\alpha$ and $\partial \Omega$.

We now define $M_o$ such that 
$$M_o\delta ^\gamma>\displaystyle\sup_{\partial \Omega_\delta\cap \Omega }u 
\ \mbox{  and } \  M_0^{1+\alpha}> {|f|_\infty\over \epsilon}.$$

By the comparison principle (Theorem \ref{comp}) 
$u^\star\leq M_od(x, \partial \Omega)^\gamma $ in $\Omega_\delta$ and the claim is
proved.

\bigskip
We now prove H{\"o}lder's regularity inside $\Omega$. 

We construct a function $\Phi$ as follows:
Let $M_o$ and $\gamma$ be as in the Claim,  
  $M= \sup  (M_o, {2\sup u\over \delta^\gamma})$ and 
 $\Phi (x) = M( |x|^\gamma)$.

We shall consider 
$$\Delta_\delta=\{(x,y)\in \Omega^2,\ |x-y|<\delta\}.$$

{\bf Claim 2} {\em For any} $(x,y)\in\Delta_\delta$
\begin{equation}\label{ho}
u^\star(x)-u_\star(y)\leq \Phi(x-y)
\end{equation}

If the Claim 2 holds this completes the proof, indeed
taking $x=y$ we would get that $u^\star=u_\star$ and then $u$ is
continuous.  Therefore going back to (\ref{ho}) 
$$u(x)-u(y)\leq {2\sup u\over \delta^\gamma} |x-y|^\gamma,$$ for
$(x,y)\in \Delta_\delta$
 which is equivalent to the local H{\"o}lder continuity. 

\bigskip

\noindent Let us check first  that (\ref{ho}) holds 
 on $\partial \Delta_\delta$. On that set,

- either $|x-y|=\delta$ and then
$u^\star(x)-u_\star(y)\leq M\delta^\gamma$ 
since  $M\delta^\gamma \geq {2\sup u}$,  

-or
$(x,y)\in\partial(\Omega \times \Omega)$. In that case, for 
$(x,y)\in (\Omega\times\partial \Omega)$ we have just proved that
$$u^\star(x)\leq M_od^\gamma\leq M|y-x|^\gamma.$$ 
Now we consider interior points. Suppose by
contradiction that $u^\star(x)-u_\star(y)> \Phi(x-y)$
for some $(x,y)\in \Delta_\delta$. Then there exists $(\bar x, \bar y)$ such that 
$$u^\star(\bar x)-u_\star(\bar y)-\Phi(\bar x-\bar y)=\sup (u^\star(x)-u_\star(y)-\Phi(x-y))>0.$$
 Clearly $\bar x\neq \bar y$.
Then using Ishii's Lemma there exists $X$ and $Y$ such that
$$( \gamma M(\bar x-\bar y)|\bar x-\bar y|^{\gamma-2}
, X)\in J^{2,+ }u^\star(\bar x)$$
$$(\gamma M(\bar x-\bar y)|\bar x-\bar y|^{\gamma-2}
, -Y)\in J^{2,-}u_\star(\bar y)$$
with 
$$\left(\begin{array}{cc}
X&0\\
0&Y
\end{array}\right) \leq \left(\begin{array}{cc}
B &-B\\
-B&B
\end{array}\right) $$
and  $B = D^2\Phi(\bar x-\bar y)$.

In particular 
$tr(X +Y)\leq 0$. We need a  more precise estimate, as in \cite{IL}. For that aim let :

$$0\leq P : = {(\bar x-\bar y\otimes \bar x-\bar y)\over |\bar x-\bar y|^2}\leq I.$$
Remarking that $X+Y\leq 4B$, one easily sees that $tr(X+Y)\leq tr(P(X+Y))\leq 4tr(PB)$.
But $tr (PB)=\gamma M(\gamma-1)| \bar x-\bar y|^{\gamma-2}<0$, hence

\begin{equation}\label{elip}
|tr (X+Y)|\geq 4\gamma M(1-\gamma)| \bar x-\bar y|^{\gamma-2}.
\end{equation}

Furthermore  by Lemma III.1 of \cite{IL} there exists a universal constant $C$ such that 
$$|X|, |Y|\leq C (|tr(X+Y)|+ |B|^{1\over 2} |tr(X+Y)|^{1\over
2}).$$
Now we can use the fact that $u$ is both a sub and a super solution of (\ref{Heq})
and applying (H2) condition 

\begin{eqnarray*}
f(\bar x)&\leq &  F(\grad_x\Phi,X)\\
&\leq &  a|\grad_x\Phi|^\alpha tr (X+Y) + 
F(\grad_y\Phi, tr(-Y))\\
&\leq & f(\bar y)    + |\grad_x\Phi|^\alpha tr (X+Y).
\end{eqnarray*}
Which implies, using (\ref{elip}), 

$$ a|\grad_x\Phi|^\alpha  4\gamma M(1-\gamma)| \bar x-\bar y|^{\gamma-2}\leq f(\bar
y)-f(\bar x). $$ Recalling that $|\grad_x\Phi|=\gamma M|\bar x-\bar
y|^{\gamma-1}$ the previous inequality becomes: 
\begin{equation}\label{M}
 aM^{\alpha+1}4\gamma^{1+\alpha}(1-\gamma)|\bar x-\bar y|^{\gamma (\alpha+1)-(\alpha+2)}\leq 2
|f|_\infty.
\end{equation} 
Using $M\geq {2(\sup u)\over \delta ^\gamma} $ and
$|\bar x-\bar y|\leq \delta$ one obtains 

$$a(2\sup u)^{1+\alpha} 4\gamma ^{1+\alpha} (1-\gamma) \delta^{-(\alpha+2)}
\leq 2|f|_\infty. $$ This is clearly false for $\delta$ small enough and it
concludes the proof.

\bigskip

\noindent {\it Proof of Theorem \ref{L}.} The proof proceeds similarly to
the proof given by Ishii and Lions in \cite{IL} but here we shall use the
fact that we already know that $u$ is H{\"o}lder continuous. 

We assume without loss of generality that in hypothesis (H2) $a= A = 1$.

Let $\mu$ be an increasing function such that $\mu(0)=0$, and
$\mu(r)\geq r$, let $l(r) = \int_0^r ds\int_0^s {\mu(\sigma)\over
\sigma} d\sigma$, let us note that since $\mu\geq 0$ for $r>0$
$$l(r)\leq rl^\prime (r)$$

Let $r_0$ be such that $l^\prime (r_0) = {1\over 2}$, $M$ such that
$Mr_0\geq 4 \sup |u|$. Let also $\delta>0 $ be given, $K= {r_0\over
\delta}$, and $z$ be such that $d(z, \partial \Omega)\geq 2\delta$. 

We define  $\varphi(x,y) =
\Phi(x-y)+L|x-z|^k$
where $\Phi(x) =M(K|x| -l(K|x|))$, 
and  
$$\Delta_z= \{ (x,y)\in \R^N\times \R^N, |x-y|<\delta, |x-z|< \delta\}.$$
We shall now choose all the constants above.

Choosing $k$ such that $k>{1\over 1-{1\over 2\gamma}}$
where $\gamma$ is such that $\gamma \in ]0,1[$ and  
$$|u(x)-u(y)|\leq c|x-y|^\gamma,$$
for some constant $c$ which depends on $u$ and $\gamma$. 

Choosing $M$ and $L$ such that 
$M\geq {2\sup u\over r_0}$
and 
$L= c\delta^{k-\gamma}$, using the H{\"o}lder continuity of $u$,
 one has  
$$u(x)-u(y)\leq \varphi(x,y)$$ on $\partial \Delta_z$. 

Suppose by contradiction that for some points $\bar x, \bar y$ one has 
$$u(\bar x)-u(\bar y)> \varphi(\bar x, \bar y).$$
Clearly $\bar x\neq \bar y$. 
Note that 
$$L|\bar x-z|^k\leq c|\bar x-\bar y|^\gamma.$$
Proceeding as in the previous proof, there exist 
$X$, $Y$ such that 

$$\left( MK(\bar x-\bar y)|\bar x-\bar y|^{-1}
(1-l^\prime (K|\bar x-\bar
y|))+ kL |\bar x-z|^{k-2} (\bar x-z), X\right)\in J^{2,+}u(\bar x)$$
and

$$\left ( MK{\bar x-\bar y\over |\bar x-\bar y|}(1-l^\prime (K|\bar x-\bar
y|)), -Y\right)\in J^{2,-}u(\bar y).$$
The matrices $X$ and $Y$ satisfy
$$\left(\begin{array}{cc}
X&0\\
0&Y
\end{array}\right) \leq \left(\begin{array}{cc}
B+\tilde L&-B\\
-B&B
\end{array}\right) $$
with $B = D^2\varphi(\bar x,\bar y)$
and

 $$\tilde L = kL|\bar x-z|^{k-2} \left(I+ (k-2){(\bar x-z\otimes \bar x-z)\over
|\bar x-z|^2}\right).$$
Let us note that  $X+Y-\tilde L\leq 4B$ and then 

$$tr(X+Y-\tilde L )\leq 4 tr(PB)$$
with 

$$P = {((\bar x-\bar y)\otimes (\bar x-\bar y))\over |\bar x-\bar y|^2}.$$
This allows to have : 

$$|tr(X+Y-\tilde L)|\geq {MK \mu(K|\bar x-\bar y|))\over |\bar x-\bar y|}\geq
MK^2.$$
Furthermore, as in the previous proof, one has 
$$|Y|\leq C (|B|^{1\over 2} |tr(X+Y-\tilde L)|^{1\over 2} + |tr(X+Y-\tilde
L)|).$$
Let us note that 
$$|B|\leq C {1\over |\bar x-\bar y|}$$
and 
  that 
$$\nabla_x \varphi(x) = MK(1-l^\prime (K|\bar x-\bar y|)){\bar x-\bar y\over
|\bar x-\bar y|} +k L|\bar x-z|^{k-2} (\bar x-z).$$
We can now use
$$L|\bar x-z|^{k-1} = (L|\bar x-z|^k)^{k-1\over k}
L^{1\over k} \leq (c\delta^\gamma)^{k-1\over k}
(c\delta^{\gamma-k})^{1\over k} = O(\delta^{\gamma-1}) =o(K)$$ so that 
for $K$ large enough 

$${2MK}\geq |\nabla_x \varphi| \geq {MK\over 4}.$$ 
Then, using the fact that $u$ is both a sub and a super solution, 
there exist some universal constants $c_1, c_2, c_3$, such that 
\begin{eqnarray*}
(MK)^\alpha|tr(X+Y-\tilde L)|&\leq &c_1(MK)^{\alpha-1} L|\bar
x-z|^{k-1} {1\over |\bar x-\bar y|^{1\over 2}}(|tr(X+Y-\tilde L)|
^{1\over 2})
+\\
&+& c_2(MK)^\alpha L|\bar x-z|^{k-1} |tr(X+Y-\tilde L)|+ f(\bar y)-f(\bar x) + \\
&+&
c_3|MK|^\alpha L|\bar x-z|^{k-2}.
\end{eqnarray*}
We shall now prove that for
 $K$ large enough this is absurd by obtaining the following estimates :

$[K1]\ L|\bar x-z|^{k-2} = O(\delta^{\gamma-2}) = o(K^2)$

$[K2]\ L|\bar x-z|^{k-1} \displaystyle{|\nabla_x\varphi|^{\alpha-1}
|tr(X+Y-\tilde L)|^{1\over 2}\over |\bar x-\bar y|^{1\over 2}} \leq o(1)
|\nabla_x
\varphi|^\alpha|tr(X+Y-\tilde L)|$

or equivalently

$[K3]\ L\displaystyle{|\bar x-z|^{k-1}\over |\bar x-\bar y|^{1\over 2}}
\leq o(1) |\nabla_x
\varphi| (K^2)^{1\over 2}= o(1) K^2.$

\bigskip

\noindent We prove $[K1]$: 
$$L|\bar x-z|^{k-2} \leq (L|\bar x-z|^k)^{k-2\over k}L^{2\over k}\leq c\delta^{\gamma-2} \leq
o(K^2).$$ 
We prove $[K3]$
\begin{eqnarray*}
{L|\bar x-z|^{k-1}\over |\bar x-y|^{1\over 2}}&\leq& L^{1\over
k}{ (L|\bar x-z|^k)^{k-1\over k}\over |\bar x-\bar y|^{1\over 2}} \\
&\leq&
CL^{1\over k} (|\bar x-\bar y|)^{\gamma (1-{1\over k})-{1\over 2}} \\
&\leq &C
(c\delta^{\gamma-k})^{1\over k} \delta^{\gamma (1-{1\over k})-{1\over 2}}
\\ &=& O(K^{{3\over 2}-\gamma} )\\
&=&o( K^{3\over 2}).
\end{eqnarray*}
We have obtained 
\begin{eqnarray*}
CK^{\alpha+2}&+&{|\nabla_x \varphi|^\alpha |tr(X+Y-\tilde L)|\over 2} \leq |\nabla_x \varphi|^\alpha |tr(X+Y-\tilde L)|\\
&\leq &2 \sup
|f|+ o(1)|\nabla_x
\varphi|^\alpha |tr(X+Y-\tilde L)|+ o(1) |\nabla_x \varphi|^\alpha +\\
& & +
o(1) K |tr(X+Y-\tilde L)|^{1\over 2} |\nabla_x\varphi|^{\alpha} \\
&\leq& |\nabla_x
\varphi|^\alpha |tr(X+Y-\tilde L)|+ o(1) K^2 |\nabla_x \varphi|^\alpha
\end{eqnarray*}
which is  a contradiction for $K$ large.
 
We have proved that for all $x$ such that 
$d(x,\partial \Omega)\geq 2\delta$ and for $y$ such that $|x-y|\leq \delta$ 
$$u(x)-u(y)\leq {2\sup M\over r_0} {|x-y|\over \delta}.$$
The local Lipschitz continuity is proved.

\section{Existence results}

\subsection{The case  $\lambda< \bar\lambda$}

In this subsection we shall prove the esistence of solutions via Perron's method by constructing explicitely a positive super solution.

\begin{theo}\label{propflambda}
Suppose that $f$ is bounded and $f\leq 0$ on $\overline{\Omega}$. 
Then, for $ \lambda<\bar\lambda$ there exists $u$  a nonnegative viscosity
solution of  
$$\left\{\begin{array}{lc}
F(\grad u,D^2u)+\lambda u^{1+\alpha}  = f &  {\rm in}\ \Omega\\
u=0\ & {\rm on}\ \partial\Omega.
\end{array}\right.
$$
Furthermore $u$ is unique.
\end{theo}

To prove this theorem, we need the two following propositions
\begin{prop}\label{proex}
Suppose that $f$ is bounded, $f\leq 0$ and $\lambda\in \R$.
Suppose that there exists $v_1\geq 0$ and $v_2\geq 0$ respectively sub solution and super solution of 
\begin{equation}\label{eqex}
\left\{\begin{array}{lc}
F(\grad v,D^2v) +\lambda v^{1+\alpha} = f & {\rm in}\ \Omega\\
v = 0 &{\rm on } \ \partial \Omega
\end{array}
\right.
\end{equation}
with $v_1\leq v_2$. Then 
 there exists a  viscosity  solution $v$ of (\ref{eqex}), such that $v_1\leq v\leq v_2$.
Moreover if $f<0$ inside $\Omega$ the solution  is unique. 
\end{prop}

\begin{prop}\label{propf}
For any  $f$ bounded  and non positive  in  $\overline{\Omega}$, there exists a unique viscosity  solution $w$  of 
\begin{equation}\label{eqf}\left\{
\begin{array}{lc}
 F(\grad w,D^2w)= f & {\rm in}\ \Omega\\
w=0 & {\rm on}\ \partial\Omega.
\end{array}\right.
\end{equation}
Of course $w$ is nonnegative by the maximum principle and H{\"o}lder continuous. 
\end{prop}

By  Proposition \ref{proex}, Proposition \ref{propf} will be proved if we  construct  a sub and super solution
for (\ref{eqf}).  
Since the null function is clearly a sub solution, it is sufficient to construct a viscosity solution $u$ of $F(\grad u,D^2u)\leq -1$
which is positive and zero on the boundary, then multiplying by the right constant we get the required super solution of (\ref{eqf}).

In the next lemma we construct such a super solution:.

\begin{lemme}
Let $\Omega$ be a bounded  ${\cal C}^2$ domain in $\R^N$. Let $d(x)=d(x, \partial \Omega)$ be the 
distance to the boundary. Then there exist $k\in \N,$ 
$\gamma\in (0,1)$, and $\beta>0$ such that 
$$u(x) = \beta\left(1-{1\over (1+d(x)^\gamma)^k}\right)$$
is a viscosity super solution of 
$$F(\grad u,D^2u)\leq -1. $$
\end{lemme}

The proof of this lemma is postponed to the appendix together with some properties of the distance function, while the proof of Proposition \ref{proex} is
 at the end of this section.

\bigskip

\noindent{\it Proof of Theorem  \ref{propflambda}}

For $\lambda<0$, one can apply directly Proposition \ref{proex}, since $0$ is a sub solution for (\ref{eqex})) and the solution constructed in Proposition
\ref{propf} is a super solution. 

\bigskip
 We now treat the case $\lambda>0$. 

We define the sequence $u_n=T_f^n(0)$ where $T_f(u)$ is defined as the unique viscosity solution of 
$$\left\{\begin{array}{lc}
F(\grad T_f(u), D^2 T_f(u)) = f-\lambda u^{1+\alpha} & {\rm in }\ \Omega\\
T_f(u)=0 & {\rm on }\; \partial \Omega.
\end{array}
\right.$$
Proposition \ref{propf} implies that $T_fu$ is well defined.

By the comparison principle and the maximum principle for $F$   in
\cite{BD},  $u_n$ is increasing and nonnegative.  We want to prove that
it is bounded. Suppose not, then $w_n := {u_n\over |u_n|_\infty}$ 
satisfies
$$F(\grad w_{n+1},D^2w_{n+1})+ \lambda \left({u_n^{1+\alpha}\over
|u_{n+1}|_\infty^{1+\alpha} }\right)= {f\over
|u_{n+1}|_\infty^{1+\alpha}}.$$ Furthermore 
$$F(\grad w_{n+1},D^2w_{n+1})+\lambda w_{n+1}^{1+\alpha} =  \lambda \left({u_{n+1}
^{1+\alpha}\over |u_{n+1}|_\infty^{1+\alpha} }-{u_n^{1+\alpha}\over
|u_{n+1}|^{1+\alpha} }\right)+ {f\over |u_{n+1}|_\infty^{1+\alpha}}\geq
{f\over |u_{n+1}|_\infty^{1+\alpha}}.$$ 
Clearly
 $$ \left |\lambda \left({u_{n+1}
^{1+\alpha}\over |u_{n+1}|_\infty^{1+\alpha} }- {u_n^{1+\alpha}\over |u_{n+1}|_\infty^{1+\alpha} }\right)+{f\over
|u_{n+1}|_\infty^{1+\alpha}}\right|\leq 2\lambda+\frac{|f|}{|u_1|_\infty^{1+\alpha}}$$
since
$0\leq {u_n^{1+\alpha}\over |u_{n+1}|_\infty^{1+\alpha} }\leq 1$.

 We can now apply  the
H{\"o}lder estimates in the previous section and this implies that the sequence $w_n$ is relatively compact 
in ${\cal C} (\bar\Omega)$.

Extracting a subsequence from $(w_n)$ and passing to the limit one gets in
particular 
$$F(\grad w,D^2w)+\lambda w^{1+\alpha} \geq 0.$$
Moreover $w=0$ on the boundary . 

We are in the hypothesis that $\lambda<\bar\lambda$ hence we can apply the maximum principle and conclude that
  $w\leq 0$. We have reached a contradiction 
since $w\geq 0$ and $|w|_\infty=1$.   

We have obtained that the sequence $u_n$ must be bounded. Since it is increasing and bounded it converges and the
convergence is uniform on  $\overline{\Omega}$, by the H{\"o}lder estimates.
Using the properties of uniform limit of viscosity solutions one gets 
that the limit
$u$ is a nonnegative  solution of 
$$F(\grad u,D^2u)+\lambda u^{1+\alpha} = f.$$

\bigskip

\noindent{\it Proof of Proposition \ref{proex} :} 
The proof relies on  Perron's method applied to viscosity solutions by Ishii
(see \cite{I}).

Let us define 

$$v = \sup\{v_1\leq u\leq v_2, u \ \mbox{is a viscosity sub solution of}\ (\ref{eqex})\}.$$
We want to prove first that $v^\star$  a sub solution.
Let $u_n$ be an increasing sequence of sub solutions, $v_1\leq u_n\leq v_2$, $u_n$ converging to $v$. 

Suppose first that  $v$ is equal to a constant $C$ on a ball $B(\bar x,r)$.  Since $C\geq 0$, it is a sub solution.

We now treat the points  where $v$ is not locally constant. Suppose by contradiction that $\bar x$ and $\varphi$ are  such
that $\nabla \varphi (\bar x)\neq 0$ and 
$$(v-\varphi)(x) \leq (v^\star-\varphi)(\bar x)=0,$$
and that there exists $r>0$ with 

$$F(\grad \varphi,D^2\varphi)(\bar x) + \lambda \varphi (\bar x)^{1+\alpha} \leq f(\bar x)-r.$$
 Let $\delta$ be small enough that for   
$|\bar x-y|\leq \delta,$  the following inequalities 
hold 
$$|F(\grad \varphi,D^2\varphi) (y)-F(\grad \varphi,D^2\varphi) (\bar x)|\leq {r\over 4},$$
$$|\varphi (y)^{1+\alpha}-\varphi (\bar x)^{1+\alpha}|\leq {r\over
4\lambda},$$
$$|f(y)-f(\bar x)|\leq {r\over 4}.$$
One can assume that the supremum of $v^\star-\varphi$ on $\bar x$ is strict, so that there exists $\alpha_\delta>0$ with  

$$\sup_{|y-\bar x|\geq \delta} (v^\star-\varphi)\leq -\alpha_\delta.$$
Finally take $N$ large enough in order that by the simple convergence of $u_n(\bar x)$ toward $v(\bar x)$ one has 
$$u_n(\bar x)-v^\star(\bar x)\geq -{\alpha_\delta\over 4 }$$
then  
$$\sup_{|x-\bar x|\leq \delta} (u_n-\varphi) (x)\geq {-\alpha_\delta \over 4} \geq -\alpha_\delta \geq  \sup_{|x-\bar x|\geq
\delta}(v^\star-\varphi)(x)\geq \sup_{|x-\bar x|\geq
\delta}(u_n-\varphi)(x).$$ Furthermore the supremum is achieved inside $B(\bar x, \delta)$, on some $x_n$. Then one has 
\begin{eqnarray*}
f(\bar x)-r&\geq& F(\grad \varphi,D^2\varphi)(\bar
x)+\lambda
\varphi(\bar x)^{1+\alpha}\\
&\geq& F(\grad \varphi,D^2\varphi) (x_n)+\lambda
\varphi (x_n)^{1+\alpha}-{r\over 2}\\
&\geq& f(x_n)-{r\over 2}\geq f(\bar
x)-{3r\over 4},
\end{eqnarray*}
  a contradiction. 
 
We now prove that $v_\star$ is a super solution. If not there would exist $\bar
x\in \Omega$, $r>0$ and $\varphi\in {\cal C}^2(B(\bar x,r)$, with  $\nabla
\varphi(\bar x)\neq 0$, satisfying 
$$0= (v_\star-\varphi)(\bar x) \leq (v_\star-\varphi)(x)$$
on $B(\bar x,r)$ such that
$$F(\grad \varphi,D^2\varphi) (\bar x) +\lambda \varphi (\bar x)^{1+\alpha} >f(\bar x).$$

We prove first that $\varphi(\bar x)< v_2(\bar x)$. If not one would have $\varphi(\bar x) = v_\star(\bar x) = v_2(\bar x)$ and then 
$$(v_2-\varphi)(x)\geq (v_\star-\varphi)(x)\geq (v_\star-\varphi)(\bar x)
= (v_2-\varphi) (\bar x)=0,$$ hence since $v_2$ is a super solution and
$\varphi$ is a test function for $v_2$ on $\bar x$,
$$F(\grad \varphi,D^2\varphi) (\bar x) +\lambda \varphi (\bar x)^{1+\alpha}\leq f(\bar x),$$
a contradiction. 
Then  $\varphi(\bar x)< v_2(\bar x)$. We construct now a sub solution which is greater than $v$ and less than $v_2$.

Let $\varepsilon>0$ be such that 
$$F(\grad \varphi,D^2\varphi)(\bar x)+\lambda \varphi(\bar x)^{1+\alpha} \geq
 f(\bar x)+\varepsilon.$$
Let $\delta$ be such that for $|x-\bar x|\leq \delta$

$$|F(\grad \varphi,D^2\varphi)(x)-F(\grad \varphi,D^2\varphi)(\bar x)|+ |f(x)-f(\bar x)|+\lambda |\varphi(x)^{1+\alpha}-\varphi(\bar x)^{1+\alpha}|\leq {\varepsilon\over 4}.$$
Then 
$$F(\grad \varphi,D^2\varphi)(x)+\lambda \varphi^{1+\alpha} 
(x) \geq f(x)+{\varepsilon\over 4}.$$
One can assume that $$(v_\star-\varphi)(x)\geq |x-\bar x|^4.$$ 
We take $r< \delta^4$ and such that $0<r< \inf_{|x-\bar x|\leq \delta}(v_2(x)-\varphi(x))$ and define 
$$w= \sup (\varphi(x)+r, v_\star)$$
$w$ is LSC as the supremum of two LSC functions.
  
One has $w(\bar x) = \varphi(\bar x)+ r$, 
and $w= v$ for $r<|x-\bar x|< \delta$.
 
$w$ is a sub solution, since when $w= \varphi+r$ one can use $\varphi+r$ as a test function,  
and since $\varphi (x)>0,$
$$F(\grad \varphi,D^2\varphi)(x)+\lambda (\varphi(x)+r)^{1+\alpha} \geq (F(\grad \varphi,D^2\varphi)+\lambda \varphi^{1+\alpha})(x) \geq f+{\varepsilon\over 4}.$$
Elsewhere $w=v$, hence it is a sub solution. Moreover $w\geq v$, $w\neq v$ and $w\leq g$. This contradicts the fact that $v$ is the
supremum of the sub solutions.  Using H{\"o}lder regularity we get that $v$ is H{\"o}lder and hence $v^\star=v_\star$.

\subsection{The case $\lambda = \bar \lambda$.}

\begin{theo}\label{eigenfunction}

Let $L$ be as in the previous section. Then, there exists $\phi>0$ in $\Omega$ such that $\phi$ is a viscosity solution of 

$$
\left\{\begin{array}{lc}
F(\grad \phi,D^2\phi)+\bar\lambda \phi^{1+\alpha} = 0 & {\rm in}\  \Omega\\
\phi=0  & {\rm on}\  \partial \Omega.
\end{array}
\right.
$$
Moreover $\phi$ is $\gamma$-H{\"o}lder continuous for all $\gamma\in ]0,1[$ 
and locally Lipschitz. 
\end{theo}
{\it Proof of Theorem \ref{eigenfunction}}

Let $\lambda_n$ be an increasing sequence which converges to $\bar\lambda$. Let $u_n$ be a nonnegative viscosity solution of $$
\left\{\begin{array}{lc}
F(\grad u_n,D^2u_n)+\lambda_n u_n^{1+\alpha} = -1 & {\rm in}\  \Omega\\
u_n=0  & {\rm on}\  \partial \Omega.
\end{array}
\right.$$

By  Theorem \ref{propflambda} the sequence $u_n$ is well defined.
We shall prove that $u_n$ is not bounded. Indeed suppose by contradiction that it is. 
Then, by  H{\"o}lder's estimate, one has that a subsequence, still
denoted $u_n$, tends uniformly to a nonnegative  continuous function $u$
which is a viscosity solution of 
$$F(\grad u,D^2u)+\bar{\lambda} u^{1+\alpha}= -1.$$
This  contradicts the definition of $\bar\lambda$. Indeed $u >0$ and one can choose 
$\varepsilon$ small enough that
$$F(\grad u,D^2u) +(\bar\lambda+\varepsilon) u^{1+\alpha}\leq -1+\varepsilon u^{1+\alpha}\leq 0.$$

We have obtained that the  sequence $|u_n|_\infty\rightarrow +\infty$. Then defining $w_n = {u_n\over |u_n|_\infty}$ one has 
$$F(\grad w_n,D^2w_n)+ \lambda_n w_n^{1+\alpha} = {f\over |u_n|^{1+\alpha}}$$
and then extracting as previously a subsequence which converges
uniformly, one gets that there exists $w$, $|w|_\infty=1$ and 
$$F(\grad w,D^2w)+\bar \lambda w^{1+\alpha} = 0.$$
The boundary condition is given by the uniform convergence.
Clearly $w$ is H{\"o}lder and locally Lipschitz continuous. 

\section{Appendix: Properties of the distance function. }

In all this section $\Omega$ is a bounded  ${\cal C}^2$ domain in $\R^N$. 
For completeness sake we shall  study the regularity of the distance function 
at points $x\in\Omega$  for which the distance to the boundary is achieved by one single point
of the boundary. See e.g. \cite{F,Fo,KP} for other interesting results on the distance function.

We shall denote the elements of
$\overline\Omega$ by $(x',x_N)\in\R^{N-1}\times\R$.  Without loss of  generality we can suppose that
$(0,0)\in \partial\Omega$ and  that $x_o= (0,d)\in \R^{N-1}\times \R^+$
is at the distance $d$ to
$\partial \Omega$ and that  $(0,0)$ is the unique point of $\partial \Omega$  at which the distance is 
achieved. 

Always without loss of generality we can suppose that there exist a neighborhood $V$ of $(0,0)\in
\R^N$ ,  $r>0$ and a function $a\in {\cal C}^2(B^\prime
(0,r))$ (the ball of center $0$ and radius $r$ in $\R^{N-1}$)   such that 
$$\partial\Omega \cap V= \{ (x^\prime, a(x^\prime)),\ x^\prime\in
B^\prime (0,r)\}.$$
And we can suppose that the unit interior normal to $\partial\Omega$ at
$(0,0)$ is
$e_N$, which implies that  
$\nabla a(0) = 0.$  
In the rest of the section we shall consider this setting.

\begin{lemme}\label{distest} 
In the above hypothesis
$$D^2 a(0)< {1\over d} I,$$
(where $I$ denotes the $N-1$ dimensional matrix identity, in the sense of
positivity of symmetric matrices.) 
 \end{lemme}

\begin{prop}\label{propdist}
In the same hypothesis there exists a neighborhood
$V_0$ of
$(0,0)$ in $\partial \Omega\times \R^+$ on a neighborhood $V_1$ of $x_o$
such that $\forall x\in V_1$ there exists one and only one $y\in V_0$
such that 
$$ |x-y| = d(x,\partial \Omega).$$
Moreover the map $x\rightarrow y(x)$ is ${\cal C}^2$ in this neighborhood. 

\end{prop}
{\it Proof of Lemma \ref{distest}}

Since $(0,0)$
is the unique point on the $N-1$ surface $x_N= a(x^\prime)$, at which the
distance is achieved, we get for
$|x^\prime|<r$ 

\begin{equation}\label{xpr}
|(0-x^\prime, d-a(x^\prime))|^2> d^2.
\end{equation}
Using 
$$a(x^\prime) = {1\over 2} \langle D^2 a (0)x^\prime, x^\prime\rangle
+ o(|x|^2),$$
(\ref{xpr}) implies that 
$$|x^\prime|^2 -d\langle D^2a(0).x^\prime,x^\prime\rangle  \geq
C|x^\prime|^4+ o(|x|^2).$$
This gives the required result.

\bigskip

\noindent{\it Proof of Proposition \ref{propdist}. }

We define a map on a neighborhood of $(0,d)$ as follows 
$$\begin{array}{rcl}
\Psi:B^\prime (0,r) \times (d-\frac{d}{2},d+\frac{d}{2})&\longrightarrow& \Omega\\
(y,t) &\longrightarrow &(y-t{\nabla a (y)\over (1+|\nabla a|^2)^{1\over 2}}\ ,\
a(y)+{t\over (1+|\nabla a|^2)^{1\over 2}}).
\end{array}
$$
 We want to prove that this map is
invertible around $(0,d)$, since  
$$D^2a (y)< (1/d)I.$$
For that aim we introduce 
$$X^\prime =y-t{\nabla a (y)\over (1+|\nabla a|^2)^{1\over 2}}$$
$$X_N =a(y)+{t\over (1+|\nabla a|^2)^{1\over 2}}$$
and we prove that the Jacobian is non zero on $(0,0)$.
A simple computation gives 
$${\partial X_i\over \partial y_j} = \delta_{ij} -t {a_{,ij}\over
(1+|\nabla a|^2)^{1\over 2} } +{ta_{,k}a_{,i}a_{,kj}\over (1+|\nabla
a|^2)^{3\over 2}}$$ for all $i,j\in [1,N-1]$, and 
$${\partial X_N\over \partial y_j} = a_{,j}+ {t (-a_{,k}a_{,kj})\over
(1+|\nabla a|^2)^{3\over 2}} $$
$${\partial X_i\over \partial t} = -{\nabla a\over (1+|\nabla
a|^2)^{1\over 2}}$$
$${\partial X_N\over \partial t} = {1\over  (1+|\nabla
a|^2)^{1\over 2}}.$$
From this one gets using $\nabla a (0) = 0$ that the Jacobian at $x_0=
(0,d)$ has the value of the determinant of the $N-1$ dimensional matrix 
$I-d (D^2a)(0)$. The previous lemma implies that this determinant is strictly positive. 
Hence in a neighborhood of $(0,d)$ the map $\psi : (y,t)\mapsto (X^\prime
, X_N)$ is invertible by the local inversion theorem.

Precisely there exists a neighborhood $V_1$ of $(0,d)$ such that for any $x\in V_1$
there exists a unique $y=y(x)\in B^\prime(0,r)$ and a unique $t$ such that:
$d(x)=t=|x-(y(x),a(y(x)))|$. Clearly $y(x)$ is differentiable and
$$\nabla y(x_0) = (I-dD^2 a)^{-1}.$$ 

This ends the proof .

\begin{cor}\label{co}
The distance $d$ is ${\cal C}^2$  around every point $x$ on which the distance is achieved in a unique point of the boundary
. Moreover in our setting at  $x_o$ the eigenvalues of $D^2d$ are
0 and the eigenvalues of $(D^2a)(I-dD^2a)^{-1}$. 
\end{cor}

{\em Proof of Corollary \ref{co}.}
We still consider the geometry and the notations of the proof of 
Proposition
\ref{propdist}. It is easy to see that for $x=(x',x_N)$ 
$$y(x)=x'- (x_N-a(y(x))\grad a(y(x)).$$
This will allow us to compute explicitly $D^2 d$. Indeed, one gets 

$$\grad y(x)=I-(x_N-a(y))D^2a.\nabla y +(\grad a\otimes \grad
a)
.\grad y.$$

Hence, in particular for $x_o=(0,d)$ with $a(0)=0$ and $\grad a(0)=0$
therefore
$\grad y(x_o)=I- dD^2a.\nabla y(x_0)$. 

Recalling that 
$$\grad d(x)=\frac{1}{d}(x'-y(x),x_n-a(y(x)))$$
then
$$D^2d(x)=\frac{1}{d}\left(\begin{array}{cc}
            I-\grad y & -\grad y.\grad a(y)\\
              O       & 1
             \end{array}\right) - \frac{1}{d^3}(x-(y,a(y))\otimes
(x-(y,a(y)))$$
where $\nabla y.\nabla a = a_{,j}y_{j,i}$. 
Then if $x_o=(0,d)$, $\grad a=0$ 
$$D^2d(x_o)=\left(\begin{array}{cc}
                      D^2a\grad y & 0 \\
                       0 &  \frac{1}{d} \end{array}\right) - \frac{1}{d^3}x_o\otimes x_o$$
where $D^2a(y) \nabla y$ is the usual product of matrices in
$\R^{N-1}$. Clearly
$D^2d(x_o)x_o=0$ so one of the eigenvalue is 0, while  for any
$x_1=(x_1^\prime,0)$ one has
$$ D^2d(x_o)x_1=\left(\begin{array}{c}D^2a\grad y x_1'\\0\end{array}\right).$$
Using the fact that $\nabla y(x_0) = (I-dD^2a(0))^{-1}$ we get the result,
choosing $x_1^\prime$ as an eigenvector of $D^2 a $. This conclude the proof.

For completeness sake let us recall that 
\begin{prop}\label{pl}
Suppose that $x\in \Omega$ is such that the distance $d$ to $\partial
\Omega$ is achieved at least on two points. Then the set
$J^{2,-}d(x)=\emptyset$.
\end{prop}

{\it Proof of Proposition \ref{pl}.}

Suppose that $x=0$ and let $y_1$ and $y_2$ be  two distinct points in
$\partial \Omega$ such that 
$d(0,\partial \Omega)= d=|0-y_1|= |0-y_2|$. 
It is sufficient to prove that $J^{2,-}d^2(0)$ is empty. 
Suppose that $a$ and $A$ are in $\R^N\times S^N$ such that for all $x$
in a neighborhood of $0$
$$d^2+a.x+ ^txAx\leq d(x, \partial \Omega)^2$$
In particular this must be satisfied for all $x = ty_1$ and $|t|< r$
small enough. This implies in particular 
$$d^2+t(a.y_1)+ t^2 (Ay_1, y_1)\leq \inf_{|t|<r} (|ty_1-y_1|^2,
|ty_1-y_2|^2)$$
In particular one gets first 
$$(a.y_1 ) t\leq -2d^2 t$$
which implies $a.y_1 = -2d^2$
and secondly one has 
$$(a.y_1)t\leq -2(y_1.y_2)t$$
which implies that $(a.y_1) = -2(y_1.y_2) = -2d^2$,
a contradiction since $y_1\neq y_2$ implies that $y_1.y_2 \neq d^2$.

\begin{prop}\label{propbeta}

Let $\Omega$ be a bounded open ${\cal C}^2$ domain in $\R^N$. Then for all
constant $\beta<0$ there exists a function $u$ which is a viscosity
solution of 
$$\left\{\begin{array}{ccc}
F(\grad u,D^2u)&\leq& \beta\ {\rm in} \ \Omega\\
u&=0&\ {\rm on}\ \partial \Omega
\end{array}
\right.
$$
$u=0$ on the boundary.
\end{prop}

{\it Proof of proposition \ref{propbeta}}
According to the previous proposition, it is enough to consider a point
$x$ where $d$ is achieved on only one point $x_o$. Hence we can consider
that the setting  is the one considered in the previous propositions.
 
Let  $K> diam \Omega$. Then $d\leq K$. Let $\gamma \in ]0,1[$ and let
$k$ be large enough to be chosen later.
% in order that 
%$${(k+2-\gamma) \gamma K^{\gamma-1}\over (1+K^\gamma
%)}2C$$
We construct the following function 
$$u(x) = 1-{1\over (1+d(x)^\gamma)^k}.$$
Clearly 
$u=0$ on the boundary, $u$ is ${\cal C}^2$ on the points where $d$ is
achieved on a unique point. 
One has 
$$\nabla u = {k\gamma d^{\gamma-1}\nabla d\over (1+d^\gamma)^{k+1}}$$
and 
$$D^2 u = {k\gamma d^{\gamma-2}\over (1+d^\gamma)^{k+2}} \left[
(\gamma-1-(k+2-\gamma) d^\gamma)\grad d\otimes\grad d + d(1 + d^{\gamma})D^2d\right].$$
We need to evaluate the eigenvalues of $D^2u$.
Using the fact that  $D^2d\grad d=0$, we obtain that
$$D^2 u.\grad d = {k\gamma d^{\gamma-2}\over (1+d^\gamma)^{k+2}} 
(\gamma-1-(k+2-\gamma) d^\gamma)\grad d,$$ hence  ${k\gamma d^{\gamma-2}\over (1+d^\gamma)^{k+2}} 
(\gamma-1-(k+2-\gamma) d^\gamma)$ is a negative eigenvalue of $D^2u$.
While for $x_i=(x_i',0)$ with $x_i'$ being an eigenvector of $D^2d$ with corresponding eigenvalue $\lambda_i$
$$D^2u.x_1=  {k\gamma d^{\gamma-1}\over (1+d^\gamma)^{k+1}}\frac{\lambda_i}{1-d\lambda_i}x_1.$$
Choosing  $\la1$ to be  the greatest eigenvalue of $D^2a(0)$, 
 we obtain that 
$${\cal M}^+_{a,A} D^2u\leq {k\gamma d^{\gamma-2}\over (1+d^\gamma)^{k+2}} \left[a(\gamma-1 -d^\gamma(k+2-\gamma))+A(N-1)\frac{\la1}{1-d\la1}d(1+d^\gamma)\right]$$
and using that $d(1+d^\gamma)\leq (1+K^\gamma)K^{1-\gamma}d^\gamma$
$${\cal M}^+_{a,A} D^2u\leq{k\gamma d^{\gamma-2}\over (1+d^\gamma)^{k+2}} \left[a(\gamma-1)- d^\gamma\left(a(k+2-\gamma)-A(N-1)\frac{\la1}{1-d\la1}(1+K^\gamma)K^{1-\gamma}\right)\right].$$
We choose  $k$ such that 

$$a(k+\gamma-2)\geq 2\left(A(N-1)\left|\frac{\la1}{1-d\la1}\right|(1+K^\gamma)K^{1-\gamma}\right).$$
Recalling that $|\grad d|=1$, that $\gamma\in(0,1)$ and  $\gamma(\alpha+1)-(\alpha+2)<0$ we have obtained that 
\begin{eqnarray*}
F(\grad u,D^2u) & \leq & |\grad u|^\alpha{\cal M}^+_{a,A} D^2u\\
 &\leq & {(k\gamma)^{\alpha +1} d^{\gamma(\alpha+1)-(\alpha+2)}\over
(1+d^\gamma)^{\alpha (k+1)+k+2}}\frac{1}{2}(\gamma-1 -d^\gamma(k+2-\gamma))\leq \beta<0.
\end{eqnarray*}
This conclude the proof.

\begin{lemme}\label{fcttest}
Let $u$ be a positive bounded function inside $\Omega$, $u=0$ on $\partial
\Omega$. Then, for all $z\in \partial \Omega$ such that $u$ is not
locally constant around $z$, there exists
$C>0$ and  $\bar x\in \Omega$, such that 
  $(2C(x-\bar x), -CI)\in J^{2,-}u(\bar x)$
\end{lemme}

{\it Proof}

Suppose that $C> {2\sup u\over (d(z,\partial \Omega))^2}$ and consider 
$$\inf_{x\in \bar\Omega} \{ u(x)+ C|z-x|^2\}$$
Let $\bar x$ be a point on which the infimum is achieved. If $\bar x\in
\partial \Omega$, this contradicts the definition of $C$. Moreover since
$u$ is not locally constant around $z$ the infimum cannot be achieved on
$z$. Then one has for all $x$
$$u(x)\geq u(\bar x)+ 2C(z-\bar x, x-\bar x)-C|x-\bar x|^2$$

\noindent Isabeau Birindelli \\
Universit{\`a} di Roma ``La Sapienza'' \\
Piazzale Aldo moro, 5\\
00185 Roma, Italy \\
e mail: isabeau@mat.uniroma1.it

\medskip
\noindent Fran\c coise Demengel\\
Universit{\'e} de Cergy Pontoise,\\
Site de Saint-Martin, 2 Avenue Adolphe Chauvin\\
95302 Cergy-Pontoise, France\\

\end{document}